\documentclass[12pt]{amsart}
\usepackage{amsmath,amscd,amsthm,amsfonts, amssymb,amsxtra}
\usepackage{epsf}
\usepackage[all]{xy}
\newtheorem{thm}{Theorem}[section]  
\newtheorem*{un-no-thm}{Theorem}
\newtheorem{lem}[thm]{Lemma}         
  
\newtheorem{bigthm}{Theorem}
\newtheorem{bigcor}[bigthm]{Corollary}

\theoremstyle{definition}

\theoremstyle{definition}
\theoremstyle{definition}

\theoremstyle{remark}
\newtheorem{rem}[thm]{Remark}

\newtheorem*{acks}{Acknowledgements}
 
\newtheorem*{out}{Outline}

\begin{document}
\title{Fiber products, Poincar\'e duality and $A_\infty$-ring spectra}
\date{\today}
\author{John R. Klein}
\address{Wayne State University, Detroit, MI 48202}
\email{klein@math.wayne.edu}
\begin{abstract} For a Poincar\'e duality space
$X^d$ and a map $X \to B$,
consider the homotopy fiber product $X \times^B X$.
If $X$ is orientable with respect to 
a multiplicative cohomology theory $E$, then, after 
suitably regrading, it is shown that the $E$-homology
of $X \times^B X$ has the structure of a graded associative algebra. 
When $X \to B$
is the diagonal map of a manifold $X$, one recovers
a result of Chas and Sullivan about the 
homology of the unbased loop space
$LX$. 
\end{abstract}
\thanks{The author is partially supported by  NSF
Grant DMS-0201695.}
\maketitle
\setlength{\parindent}{15pt}
\setlength{\parskip}{1pt plus 0pt minus 1pt}
\def\Top{\bold T\bold o \bold p}
\def\Sp{\bold S\bold p}
\def\vo{\varOmega}
\def\vs{\varSigma}
\def\smsh{\wedge}
\def\flush{\flushpar}
\def\id{\text{id}}
\def\dbslash{/\!\! /}
\def\codim{\text{\rm codim\,}}
\def\:{\colon}
\def\holim{\text{\rm holim\,}}
\def\hocolim{\text{\rm hocolim\,}}
\def\hodim{\text{\rm hodim\,}}
\def\hocodim{\text{hocodim\,}}
\def\Bbb{\mathbb}
\def\bold{\mathbf}
\def\Aut{\text{\rm Aut}}
\def\cal{\mathcal}
\def\frak{\mathfrak}

\section{Introduction}

Let $M^d$ be an orientable, closed manifold. 
Using  intersection theory on singular chains,
M.\ Chas and D.\ Sullivan \cite{chas-sullivan-string} have constructed
an operation 
$$
\bullet\:H_p(LM) \otimes H_q(LM) \to H_{p+q-d}(LM)\, ,
$$
called the {\it loop product},
where $LM = \text{map}(S^1,M)$ denotes the free loop space
of $M$. If we set $H_*(LM)[d] := H_{*+d}(LM)$,
then the loop product gives $H_*(LM)[d]$ the structure of an
graded associative ring.
The action of the circle 
on $LM$ given by rotating loops gives rise to  a differential $\Delta\:
H_*(LM)[d]\to H_{*+1}(LM)[d]$. Chas and Sullivan prove that
the pair $(\bullet,\Delta)$ gives $H_*(LM)[d]$ the structure
of a {\it Batalin-Vilkovisky} algebra, i.e., a 
graded commutative algebra with differential 
such that the Leibniz rule fails up to a term which is bilinear. 
Recently, R.\ Cohen and J.\ Jones 
\cite{cohen-jones-string} gave a spectrum level description
of the loop product in terms of a Pontryagin construction. 
Shortly thereafter, using equivariant Spanier-Whitehead duality, 
W.\ Dwyer and H.\ Miller and  
the author  (independently, unpublished) described the loop product via
the identification of $H_*(LM)$ with topological Hochschild homology.
A spin-off of the last approach is that 
the loop product exists in the more general setting 
when $M$ is an orientable Poincar\'e duality space.
\medskip
 
This paper has two goals. The first is to
exhibit algebra structures on the homology
of a wider class of spaces. To each space in the class
we will associate a certain Thom spectrum.
The algebra structures
then follow by showing that the
Thom spectra are
$A_\infty$-ring spectra.
The second goal of this paper is to illustrate how some of
the already known examples (based loop spaces, orientable closed manifolds
and free loop spaces) fit into the wider class.
\medskip

To describe our class,
fix a Poincar\'e duality space $X$ of formal dimension $d$.
Let $X \to B$ be a map of based spaces, which
for the purposes of exposition we take to be a Serre fibration. Then
we consider the fiber product 
$$
X \times^B X \, .
$$ 

Let 
$-\tau$ denote the Spivak normal fibration of $X$.
This is a stable spherical fibration over $X$. By desuspending,
our convention will
be that the fiber of $-\tau$
is a $(-d)$-sphere spectrum. Let
$$
(X \times^B X)^{-\tau}
$$
be the effect of taking the Thom spectrum of 
the pullback of $-\tau$ along the first factor projection
$X\times^B X \to X$.

The main theorem of this paper is

\begin{bigthm} \label{main} Assume $X$ is
connected. Then
$$
(X \times^B X)^{-\tau}
$$
is an $A_\infty$-ring spectrum.

Consequently, if $E$ denotes a multiplicative
cohomology theory for which $X$ is $E$-orientable,
then 
$$
E_*(X\times^B X)[d] := E_{*+d}(X\times^B X)
$$
has the structure of a graded associative ring.
\end{bigthm}

(The  assumption that $X$ is connected 
is  an artifact of our method of proof: 
Theorem  \ref{main} is indeed true when $X$ is disconnected, but
we will not address this issue in the paper.)
\medskip

Theorem \ref{main} is actually a corollary of our next result which
identifies the above Thom spectrum with
a certain ring spectrum of endomorphisms. 
To formulate the result, we give $X$ a basepoint.
The map $X \to B$ then becomes a map of based spaces.
Let $F$ denote its homotopy fiber. Then $F$
may be equipped with an action by $\Omega B$, where the
latter is a suitable topological group version of the 
based loop space of $B$ (the Borel construction
of this action recovers the map $X\to B$ up to homotopy).
Let $F_+$ be the disjoint union of $F$ with a basepoint. Then the suspension
spectrum $\Sigma^\infty F_+$ comes equipped with an (naive) $\Omega B$-action. 
\bigskip

\begin{bigthm} \label{endo} Assume $X$ is connected. 
Then there is a weak homotopy equivalence of spectra
$$
(X \times^B X)^{-\tau} \,\, \simeq \,\, 
\text{\rm end}_{\Omega B}(\Sigma^\infty F_+) \, ,
$$
where the right side denotes the endomorphism 
spectrum of $\Omega B$-equivariant
stable self maps of $F_+$.
\end{bigthm}

\subsection*{Examples}
\medskip

{\flushleft (1).} {\sl $X \to B$ is the diagonal map of $X$.}
\smallskip

Then $(X\times^B X)^{-\tau}$ is weak equivalent to
$(LX)^{-\tau}$, where the latter is formed by
taking the Thom spectrum with respect to pulling back
the Spivak fibration of $X$ along the evaluation
map $LX \to X$. 
The endomorphism spectrum appearing in the statement of
Theorem \ref{endo} is the {\it topological Hochschild cohomology} 
of the $A_\infty$-ring $S[\Omega X] := 
\Sigma^\infty (\Omega X_+)$. 
In their solution to Deligne's
Hochschild cohomology conjecture, 
J.\ McClure and J.\ Smith have shown that the  
topological Hochschild cohomology 
of an $A_\infty$-ring always admits an action by an operad 
which is weak equivalent to the
little 2-disks operad \cite{MS}. 
Consequently, Theorem \ref{endo} implies that
$(LX)^{-\tau}$ admits an action
by an operad which is weak equivalent to the little  
$2$-disks. In particular, 
$(LX)^{-\tau}$ is a homotopy commutative
ring spectrum (an observation first made by Cohen
and Jones in the manifold case).
\bigskip

{\flushleft (2).} {\sl $X = B$ and the  map $X\to B$ is the identity.}
\smallskip

Since $X$ is $E$-orientable, we have by Poincar\'e duality
$$
E_*(X)[d] \,\, \cong \,\, E^{-*}(X)\, ,
$$
and this isomorphism defines the multiplication
on $E_*(X)[d]$.

In particular, when $E_*$ is singular homology,
the multiplication is given by the intersection pairing
of $X$.
\bigskip

{\flushleft (3).} {\sl $X$ is a point.}
\smallskip 

In this instance,
 $X\times^B X$ is identified with the based loop space $\Omega B$.
Theorem \ref{main} says that $E_*(\Omega B)$ is a graded ring for
any multiplicative theory $E$. The ring structure in this case is
given by Pontryagin product, i.e., the homomorphism induced by loop multiplication.
\bigskip

{\flushleft (4).} {\sl $B$ is a point}. 
\smallskip

If $E$ is singular homology and $X$ is an orientable manifold, then
the ring structure on $H_*(X\times X)[d]$ can be described on 
the chain level  as follows: 
if $\sigma \otimes \tau \in C_p(X) \otimes C_q(X)$
is a cycle 
and $\sigma'\otimes \tau' \in  C_{p'}(X) \otimes C_{q'}(X)$
is another cycle
then 
$$
(\sigma \otimes \tau)\bullet (\sigma' \otimes \tau')
\,\, := \,\, \epsilon(\tau,\sigma')\sigma \otimes \tau'
\in C_p(X) \otimes C_{q'}(X)\, ,
$$
where $\epsilon(\tau,\sigma')$ is trivial unless
$q + p' = n$, in which case $\epsilon(\tau,\sigma')$
denotes the intersection number $\tau \cdot \sigma'$. 
\bigskip

{\flushleft (5).} {\sl  $B = \text{\rm map}(K,X)$ for a CW complex
$K$ and 
$X \to \text{\rm map}(K,X)$ is the inclusion of the constant maps.}
\smallskip

Then $X \times^B X$ is identified with 
$$
\text{ map}(SK,X)\, ,
$$
where $SK$ denotes the unreduced suspension of $K$. Consequently,
$H_*(\text{map}(SK,X))[d]$
has a ring structure.

The topological monoid $\text{Aut}(SK)$ of self homotopy equivalences
of $SK$ acts by composition on $\text{map}(SK,X)$. 
On singular homology, it
would be interesting to understand how the action 
intertwines with the ring structure, since in the special
case when $K=S^0$, one knows that this information
encodes the Batalin-Vilkovisky structure on $H_*(LX)[d]$.  
\bigskip

The next result of this paper concerns
the compatibility of the $A_\infty$-ring
structures with respect to base change.
It is clear that a map of spaces $B \to C$ gives rise to
a map of Thom spectra
$$
(X\times^B X)^{-\tau} \to (X\times^C X)^{-\tau}
$$
where $X \to C$ is given by the composition
$X \to B \to C$.

\begin{bigthm} \label{compatibility} 
Assume $X$ is connected. Then the map 
$$
(X\times^B X)^{-\tau} \to (X\times^C X)^{-\tau}
$$
is a morphism of $A_\infty$-rings.
\end{bigthm}

In the special case when $B \to C$ is 
given by the projection of $X\times X$ onto its first factor, 
Theorem \ref{compatibility} reduces to an observation
already made by Cohen and Jones 
\cite[Th.\ 1(1)]{cohen-jones-string}:
the map $(LX)^{-\tau} \to X^{-\tau}$ 
induced by loop evaluation is a morphism of $A_\infty$-rings.  
\medskip

Our final result concerns Thom spectra of the
form $(X \times^B X)^\xi$ where $\xi$ is
a spherical fibration on $X \times^B X$ 
is induced from one on $X$ via base change along 
the first factor projection $X \times^B X \to X$.
As above, we assume that $X$ is a connected
Poincar\'e duality space.

\begin{bigthm}\label{module} 
The Thom spectrum 
 $$(X \times^B X)^\xi$$ 
is a left $(X \times^B X)^{-\tau}$-module. 
\end{bigthm}
\medskip

As a special case, we obtain what seems to be a new result
about the free loop space:

\begin{bigcor} $\Sigma^\infty LX_+$ is a left $(LX)^{-\tau}$-module.
\end{bigcor}
\medskip

\begin{out}
\S2 is language. In \S3 we prove Theorem \ref{endo};
the  main idea is to use the {\it norm map} 
for equivariant spectra
constructed by the author in \cite{Klein6}.
In \S4 we use Theorem \ref{endo} together with results
of \cite{EKMM} to prove Theorem \ref{main}. In \S5 we prove 
Theorem \ref{compatibility}. In \S6 we use ideas
similar \S4 and \S5 to prove Theorem \ref{module}.
\end{out}

\begin{acks} I wish to thank to Ralph Cohen, John Jones, 
Mike Mandell and Jim McClure for discussions having a bearing
on this paper. I am especially indebted to Jack Morava who 
first introduced me to the Cohen and Jones paper \cite{cohen-jones-string}.
\end{acks}

\section{Language}

This section is not intended to be complete. A more detailed
exposition of this material appears in \cite{Klein6}.
 
\subsection*{Spaces}
All spaces will be compactly generated, and
we make the convention that
products are to be re-topologized using the compactly
generated topology. Mapping spaces are to be given
the compactly generated, compact open topology.

A {\it weak equivalence} of
spaces denotes (a chain of) weak homotopy equivalence(s).

If $f\:A\to C$ and $g\:B \to C$ are maps of spaces, then
the {\it homotopy fiber product} (or {\it homotopy pullback})
is the space $A\times^C B$ consisting of triples
$(a,\lambda,b)$ with $a\in A$, $b\in B$ and $\lambda\:[0,1]\to C$
satisfying $f(a) = \lambda(0)$ and $g(b) = \lambda(1)$. If either 
$f$ or $g$ is a fibration, then the evident map from the fiber
product into the homotopy fiber product is a weak equivalence.

\subsection*{Poincar\'e spaces} 
A finitely dominated space $X$ is said to be an
 {\it orientable Poincar\'e duality space} 
of (formal) dimension $n$ if there exists
a fundamental class $[X] \in H_n(X;\Bbb Z)$
such that the associated cap product homomorphism
$$
\cap [X]\:H^*(X) \to H_{n{-}*}(X)
$$
is an isomorphism in all degrees. 
Similarly, one has the notion of Poincar\'e duality
with respect to a multiplicative cohomology theory $E$, where
a fundamental class is required to live in the abelian\
group $E_n(X) := \pi_n(E\smsh X_+)$.

\subsection*{Naive $G$-Spectra} 
Let $G$ be the geometric realization of a simplicial group.
A {\it (naive) $G$-spectrum} consists of based
(left) $G$-spaces $E_i$ for $i \ge 0$, and  
equivariant based maps $\Sigma E_i \to E_{i{+}1}$
(where we let $G$ act  trivially on the suspension coordinate of $\Sigma E_i$).
A {\it morphism} $E \to E'$ of $G$-spectra consists of
maps of based spaces $E_i \to E'_i$ which are compatible with
the structure maps.  A {\it weak equivalence}  of
$G$-spectra is a map inducing an isomorphism on
homotopy groups. $E$ is an {\it $\Omega$-spectrum} if
the adjoint maps $E_i \to \Omega E_{i{+}1}$ are weak
homotopy equivalences.

If $X$ is a based $G$-space, then its 
{\it suspension spectrum} $\Sigma^{\infty} X$ is a $G$-spectrum with
$j$-th space $Q(S^j \smsh X)$,
where $Q = \Omega^{\infty}\Sigma^{\infty}$ is the stable
homotopy functor (here $G$ acts trivially on 
the suspension coordinates). We use the notation 
$S[G]$ for the suspension spectrum of $G_+$.
considered as a  $(G{\times}G)$-spectrum (the action on
$G_+$ is given left multiplication with respect to the
first factor of $G{\times} G$ and  right multiplication composed with the
involution $g \mapsto g^{-1}$ on the second factor.

We now give some constructions on $G$-spectra. The extent 
to which each construction is homotopy invariant is indicated 
in parenthetically.

\bigskip

If $U$ is a based $G$-space and $E$ is a $G$-spectrum, then 
$$
U \smsh E
$$
denotes the $G$-spectrum which in degree $j$ is the smash
product $U \smsh E_j$ provided with the diagonal action
(this has the correct homotopy type when
$U$, considered unequivariantly, is a based CW complex). 
The associated {\it orbit spectrum} 
$$
U\smsh_G E
$$
is given by taking $G$-orbits degreewise (it
has the correct homotopy type when $U$ is a based $G$-CW complex which is
free away from the basepoint).
\bigskip

Similarly, we can form the function spectrum
$$
\text{\rm map}(U,E)
$$
which in degree $j$ is given by
$\text{\rm map}(U,E_j) =$ the 
function space of unequivariant based maps from $U$ to $E_j$. The
action of $G$ on $\text{\rm map}(U,E)$
is given by conjugation, i.e., $(g*f)(u) = gf(g^{-1}u)$ for
$g\in G$ and $f \in \text{\rm map}(U,E_j)$
(the spectrum $\text{\rm map}(U,E)$ has the correct homotopy type
when $E$ is an $\Omega$-spectrum and $U$ is a CW complex).

Let 
$$
\text{\rm map}_G(U,E)
$$
denote the {\it fixed point spectrum} 
of $G$ acting on $\text{\rm map}(U,E)$ , i.e., the
spectrum whose $j$-th space consists of the equivariant
functions from $U$ to $E_j$ (this has the correct homotopy
type if $E$ is 
an $\Omega$-spectrum and
$U$ is the retract of a based $G$-CW complex which
is free away from the base point).
\bigskip

If $E$ is a  $G$-spectrum then
the {\it homotopy orbit spectrum} $E_{hG}$
is 
$$
 E \smsh_G EG_+\, ,
$$ 
where $EG$ the free contractible $G$-space (arising from the bar construction),
and $EG_+$ is the effect of adding a basepoint to $EG$.

The {\it homotopy fixed point spectrum} $E^{hG}$ is 
$$
\text{\rm map}_G(EG_+,E) \, .
$$
\medskip

In the above constructions, the hypotheses granting the correct
homotopy type can always be achieved by changing the input
spectra up to natural weak equivalence. This follows from
a specific choice of Quillen model structure on 
the category of $G$-spectra (for details,
see \cite{Schwede}).

\subsection*{$A_\infty$-rings} Roughly,
 A {\it (strict) $A_\infty$-ring} consists of a spectrum $R$, a product
$R\smsh R \to R$ and a unit $S^0 \to E$ such that the 
axioms for a classical ring (associativity, etc.) are relaxed to 
hold up to homotopy and all higher homotopy coherences
(i.e., $R$ is an algebra over
the associahedron operad). 
More generally, we use the term $A_\infty$-ring for any spectrum
having the weak homotopy type of a strict $A_\infty$-ring.

Essentially, the only fact 
we use in this paper about $A_\infty$-rings
is that taking the (enriched) endomorphisms of an 
object in a suitably nice category of naive  $G$-spectra 
forms an $A_\infty$-ring. This will follow 
from the existence of a good smash product construction (\cite{EKMM}).

\section{The proof of Theorem \ref{endo}} 

\subsection*{Some technical simplifications} 
Once a basepoint for $X$ is chosen, the
map $X \to B$ becomes a based map. Let $B_0$ denote
the connected component of $B$ containing the basepoint.
Using the assumption that $X$ is connected,
it is straightforward to check that the map of homotopy fiber
products
$$
X\times^{B_0} X  \to  X \times^B X
$$
is a weak homotopy equivalence (we omit the argument). 
Hence the associated map of Thom spectra
is also a weak equivalence.  {\it We therefore assume henceforth 
that $B$ is a connected space.}

Let $G = \Omega B$ denote the geometric realization of
the total singular complex of the Kan loop group of $B$
(see \cite{wald-loop}; in what follows, we abbreviate terminology and call
$\Omega B$ the {\it loop group} of $B$). 
Then we have a natural weak
equivalence of based spaces $BG \simeq B$. Using this identification,
we will assume that $B$ has been replaced by $BG$. We
therefore have a map $X \to BG$. The homotopy  fiber of this map
is then identified with the (strict) fiber product
$$
F\,\, :=  \,\, EG \times^{BG} X \, .
$$ 
This description of the homotopy fiber equips it with a preferred action
of $G$ (arising from the action of $G$ on the first factor $EG$). The Borel
construction $EG \times_G F \to BG$ is then identified with the
map $X \to B$ up to fiberwise weak equivalence. 

Let $H = \Omega X$ be loop group of $X$. Then the Borel construction of
$H$ acting on $F$ gives a fibration
$$
EH \times_H F \to BH \, ,
$$
which, with respect to the identification $BH \simeq X$,
is fiberwise weak homotopy equivalent to 
the first factor projection $X \times^B X \to X$. 

We next give a Thom spectrum version of the above.
Let $S^{-\tau}$ denote the fiber of the Spivak normal
fibration of $X \simeq BH$ desuspended to down to degree $-d$. 
As above, we equip 
$S^{-\tau}$ with an $H$-action in such a way that
the Borel construction of the action coincides with the
Spivak fibration of $X$. 
Then we have a weak equivalence
of spectra
\begin{equation}\label{thomify}
S^{-\tau} \smsh_{hH} F_+ \,\, \simeq \,\, (X\times^B X)^{-\tau}
\end{equation}
where the left hand side denotes the homotopy orbit
spectrum of $H$ acting diagonally on 
$S^{-\tau} \smsh  F_+$.

\subsection*{The norm map} Let $E$ be a $G$-spectrum.
In \cite{Klein6}, the author constructed a (weak) map of spectra
$$
\eta\:D_G\smsh_{hG} E \to E^{hG}
$$
which is a natural in $E$. Here $D_G = S[G]^{hG}$ is the {\it dualizing
spectrum} of $G$ which is given by taking the homotopy fixed points
of $G$ acting by left multiplication on $S[G] =$ the suspension
spectrum of $G_+$. Right multiplication 
by $G$ composed with
the involution $g \mapsto g^{-1}$ gives $D_G$ the structure
of a $G$-spectrum.

The author also proved in \cite{Klein6} 
that $\eta$ is a weak equivalence whenever
$BG$ is a finitely dominated space. Furthermore, when $BG$ is
finitely dominated it was shown that $BG$ is a Poincar\'e duality space
of formal dimension $d$ if and only if $D_G$ is unequivariantly weak
equivalent to $S^{-d}$. In this instance, it was also shown that
$D_G$ gives a model for the {\it Spivak fiber} of $BG$, i.e., the
unreduced Borel construction of $G$ acting on $D_G$ gives 
a stable spherical fibration $EG \times_G D_G \to BG$ which
is the Spivak normal fibration of the Poincar\'e space $BG$.
In the sequel, we fix the notation $$D_G = S^{-\tau}$$ whenever $BG$
is a Poincar\'e duality space.
\bigskip 

We apply the
norm map $\eta$ to (the suspension spectrum of) $F_+$ and
the group $H$. Since $BH = X$ is a Poincar\'e duality
space, we see that $\eta$  takes the form of a weak equivalence
\begin{equation}\label{norm}
S^{-\tau} \smsh_{hH} F_+ \overset \sim \to (\Sigma^\infty F)^{hH}\, .
\end{equation}

\subsection*{Change of groups} In the following,
let $H \to G$ be the homorphism of loop groups
arising from a map of based spaces.

\begin{lem}[``Shapiro's Lemma'']
\label{change} Let $W$ be a $G$-spectrum which is also an
$\Omega$-spectrum. Let $Y$ be a based $H$-CW complex whose action
is free away from the basepoint.

Then there is a natural weak equivalence of spectra
$$
\text{\rm map}_H(Y,W)\,\,  \simeq \,\, 
\text{\rm map}_G(Y\smsh_H G_+,W)\, .
$$
\end{lem}

\begin{proof} In fact, we claim the two spectra
are isomorphic. Taking adjoints, we get
$$
\text{\rm map}_G(Y\smsh_H G_+,W) \cong
\text{\rm map}_H(Y, \text{\rm map}_G(G_+,W))
$$
On the other hand, $W \cong\text{\rm map}_G(G_+,W)$.
\end{proof}
\medskip

We now apply \ref{change} to the homorphism arising
from the map $X \to B$,
with $Y = EH_+$ and $W = \Sigma^\infty F_+$.
The result yields a  weak equivalence of spectra
$$
(\Sigma^\infty F_+)^{hH} \,\, \simeq \,\, 
\text{\rm map}_G(EH\smsh_H G_+, \Sigma^\infty F_+) \, .
$$
Clearly, there is a weak equivalence
of based $G$-spaces $EH\smsh_H G_+ \simeq F_+$ and thus the right
side is identified with $\text{\rm map}_G(F_+, \Sigma^\infty F_+)$.
Regarding the latter as the spectrum
of equivariant stable self-maps of $F_+$, 
let us substitute the notation
$$
\text{\rm end}_{G}(\Sigma^\infty F_+) \,\, := 
\,\, 
\text{\rm map}_G(F_+, \Sigma^\infty F_+)\, .
$$
Then \eqref{thomify} and \eqref{norm} 
above yield a weak equivalence of spectra
$$
(X\times^B X)^{-\tau} \,\, \simeq \,\, \text{\rm end}_G(\Sigma^\infty F_+)\, .
$$
This completes the proof of Theorem \ref{endo}.

\section{The proof of Theorem \ref{main}}

By Theorem \ref{endo}, it is enough to check
that $\text{\rm end}_G(\Sigma^\infty F_+)$ is an $A_\infty$-ring.
We will simply quote \cite{EKMM} to prove this.

Recall that $S[G] := \Sigma^\infty (G_+)$. Then $S[G]$ is
an $S$-algebra and $\Sigma^\infty F_+$ is a left $S[G]$-module
(see \cite[IV.\ Th.\ 7.8]{EKMM}). 
Furthermore, there is an
evident identification 
$$
\text{\rm end}_G(\Sigma^\infty F_+) 
\simeq 
\text{\rm end}_{S[G]}(\Sigma^\infty F_+) \, ,
$$
where the right side denotes the function object
of $S[G]$-module self maps of $\Sigma^\infty F_+$.
By \cite[III.\ Prop.\ 6.12]{EKMM}, 
$\text{\rm end}_{S[G]}(\Sigma^\infty F_+)$ is an
$S$-algebra.  Finally, by \cite[II.\ Lem.\ 3.4]{EKMM},
one knows that an $S$-algebra is an $A_\infty$-ring.
\medskip

\begin{rem}  A nuts and bolts argument can be given
to prove the weaker statement that
$E := \text{\rm map}_G(F_+,\Sigma^\infty F_+)$ is a ring spectrum
(i.e., without verifying the $A_\infty$ condition).
We now sketch this argument.

The $j$-th space of $E$ is
$$
E_j = \text{\rm map}_G(F_+,Q(\Sigma^jF_+)) \, ,
$$
where $Q = \Omega^\infty \Sigma^\infty$ denotes 
the stable homotopy functor.

For nonnegative integers $i$ and $j$, there is a map
$$
E_i \smsh E_j \to 
E_{i+j}
$$
which is given as follows: by taking adjoints, a point
in $E_i$   amounts to 
a map of spectra $\Sigma^\infty F_+ \to \Sigma^\infty\Sigma^i F_+$.
Using in also the evident homeomorphism
$Q(\Sigma^jF_+) \cong \Omega^i Q(\Sigma^{i{+}j}F_+)$, 
a point in $E_j$  amounts to a map
of spectra $\Sigma^\infty \Sigma^i F_+ \to \Sigma^\infty \Sigma^{i{+}j}F_+$.
The composition of these maps of spectra is then a map
$\Sigma^\infty F_+ \to  \Sigma^\infty \Sigma^{i{+}j}F_+$, which
is adjoint to a point of $E_{i{+}j}$. 
The maps  $E_i \smsh E_j \to E_{i{+}j}$ 
assemble to a morphism of {\it bispectra} 
$E\otimes E \to E$, where $E\otimes E$ denotes
the external smash product of $E$ with itself, and the codomain
$E$ is regarded as a bispectrum in the obvious way. This map of bispectra
induces a map of spectra
$$
E\smsh E \to E\, ,
$$
which makes $E$ into a ring spectrum.
\end{rem}

\section{The proof of Theorem \ref{compatibility}}

We rename $F = F_B$ to indicate that
it is a $\Omega B$-equivariant model for the homotopy fiber
of $X\to B$ and similarly, we have $F_C =$
an equivariant model for the homotopy fiber
of $X \to C$. Then there is an $\Omega C$-equivariant
weak equivalence of spaces
$$
F_C \,\,  \simeq  \,\, F_B \times_{\Omega B} \Omega C \, .
$$

Using Theorem \ref{endo}, the
map $(X\times^B X)^{-\tau} \to (X\times^C X)^{-\tau}$
corresponds to 
\begin{equation}\label{end-map}
\text{\rm end}_{\Omega B}(\Sigma^\infty (F_B)_+)
\to \text{\rm end}_{\Omega C}(\Sigma^\infty (F_C)_+) \, .
\end{equation}
A straightforward checking of the details of the 
proof of Theorem \ref{endo} which we omit shows that this map
is described up to homotopy by {\it induction}
along the homomorphism
$\Omega B \to \Omega C$. Explicitly, a self map
 $f\:\Sigma^\infty (F_B)_+ \to \Sigma^\infty (F_B)_+$ induces
a self map 
$$
f\smsh \text{id}_{(\Omega C)_+}\:
\Sigma^\infty (F_B)_+ \smsh_{\Omega B} (\Omega C)_+ \,\, \to \,\,
 \Sigma^\infty (F_B)_+ \smsh_{\Omega B} (\Omega C)_+ \, ,
$$
and using the identification of
$F_C$ above, we get an $\Omega C$-equivariant weak 
equivalence
$$
\Sigma^\infty (F_C)_+ \,\, \simeq \,\,
\Sigma^\infty (F_B)_+ \smsh_{\Omega B} (\Omega C)_+ \, .
$$

Hence, on the level of $S$-algebras, we can rewrite the map
\eqref{end-map} as
$$
\text{\rm end}_{S[\Omega B]}(\Sigma^\infty (F_B)_+)
\to \text{\rm end}_{S[\Omega C]}(\Sigma^\infty (F_B)_+ \smsh_{S[\Omega B]}
S[\Omega C]) \, .
$$
With $R = S[\Omega B]$ and $R' = S[\Omega C]$,
the above amounts to the question of whether 
the extension of scalars functor  $M \mapsto M \smsh_R R'$ is enriched
over $S$-modules. As pointed out to me by M.\ Mandell,
the latter 
is a formal consequence of the result that $S$-modules form a 
closed symmetric monoidal category \cite[II.\ Th.\ 1.6]{EKMM}.

\section{The proof of Theorem \ref{module}}

We argue as in the proof of Theorem \ref{endo}, and use the notation
of that proof.

Let $S^{\xi + \tau}$ denote the fiber of spherical fibration over
$BH \simeq X$ which classifies $\xi + \tau$, where $\tau$ is 
the Spivak tangent bundle of $BH$ (whose fiber is a sphere
spectrum of dimension $d = \dim X$). Then $S^{\xi + \tau}$ is
a sphere spectrum with $H$-action whose unreduced Borel construction
represents $\xi +\tau$. Hence, we have
$$
\begin{aligned}
(X{\times}_B X)^\xi \quad &\simeq \quad
S^{\xi}\smsh_{hH} F_+ ,  \quad & \text{cf.\ \S3}, \\
 \quad &\simeq \quad 
S^{-\tau}\smsh_{hH} S^{\xi +\tau}\smsh F_+  \quad & \text{by rewriting},\\
 \quad &\simeq \quad
 (S^{\xi+\tau}\smsh F_+)^{hH} \quad & \text{by } \eta \text{ (cf.\ \S3)}, \\
  \quad &\simeq \quad
\text{\rm map}_G(F_+,S^{\xi +\tau}\smsh F_+)\, \quad &\text{by change of groups.} 
\end{aligned}
$$
\bigskip

If we now translate the last term into category of $S$-modules,
it becomes identified with
$$
\text{\rm map}_{S[G]}(\Sigma^\infty F_+, S^{\xi+\tau}\smsh F_+) \, .
$$
The closed symmetric monoidal structure of $S$-modules then gives
a composition pairing
$$
\text{\rm end}_{S[G]}(\Sigma^\infty F_+) \,\, \smsh_S \,\,
\text{\rm map}_{S[G]}(\Sigma^\infty F_+, S^{\xi+\tau}\smsh F_+)
\,\, \to \,\, \text{\rm map}_{S[G]}(\Sigma^\infty F_+, S^{\xi+\tau}\smsh F_+)
$$
and we infer that    
$\text{\rm map}_{S[G]}(\Sigma^\infty F_+, S^{\xi+\tau}\smsh F_+)$ is
an $\text{\rm end}_{S[G]}(\Sigma^\infty F_+)$-module.
It follows that $(X\times^B X)^{\xi}$ is a left
$(X\times^B X)^{-\tau}$-module.





\begin{thebibliography}{E-K-M-M}

\bibitem[C-J]{cohen-jones-string}%
Cohen, R.~L., Jones, J.~D.~S.: {A homotopy theoretic realization of string
  topology}.
\newblock {\it Math.\ Annalen} {\bf 324}, 773--798 (2002)

\bibitem[C-S]{chas-sullivan-string}%
Chas, M., Sullivan, D.: String topology.
\newblock {\it MathArXiv preprint math.GT/0212358}

\bibitem[E-K-M-M]{EKMM}%
Elmendorf, A.~D., Kriz, I., Mandell, M.~A., May, J.~P.: {Rings, Modules, and
  Algebras in Stable Homotopy Theory}.
\newblock (Mathematical Surveys and Monographs, Vol.~47).
\newblock Amer.\ Math.\ Soc. 1997

\bibitem[Kl]{Klein6}%
Klein, J.~R.: The dualizing spectrum of a topological group.
\newblock {\it Math. Annalen} , 421--456 (2001)

\bibitem[M-S]{MS}%
McClure, J.~E., Smith, J.~H.: {A solution of Deligne's Hochschild cohomology
  conjecture}.
\newblock In: {Recent progress in homotopy theory (Baltimore, MD, 2000)},
  pp.~153--193.
\newblock Amer. Math. Soc. 2002

\bibitem[Sc]{Schwede}%
Schwede, S.: {Spectra in model categories and applications to the algebraic
  cotangent complex}.
\newblock {\it J. Pure Appl. Algebra\rm } {\bf 120}, 77--104 (1997)

\bibitem[Wa]{wald-loop}%
Waldhausen, F.: {On the construction of the Kan loop group}.
\newblock {\it Doc. Math.} {\bf 1}, 121--126 (1996)
\end{thebibliography}
\end{document}